\documentclass[12pt,twoside]{article}
\usepackage[english]{babel}
\usepackage[latin1]{inputenc}
\usepackage{amsmath}
\usepackage{amssymb,amsfonts}
\usepackage{graphicx}
\usepackage{times,amssymb,amscd}

\newtheorem{theorem}{Theorem}[section]

\newtheorem{exmple}[theorem]{Example}
\newtheorem{defn}[theorem]{Definition}
\newtheorem{rmrk}[theorem]{Remark}

\newcommand{\bG}{\mathbf{G}}
\newcommand{\bH}{\mathbf{H}}

\newcommand{\bN}{\mathbf{N}}
\newcommand{\bR}{\mathbf{R}}
\newcommand{\bS}{\mathbf{S}}
\newcommand{\bV}{\mathbf{V}}

\newcommand{\bs}{\mathbf{s}}
\newcommand{\ba}{\mathbf{a}}

\newcommand{\bx}{\mathbf{x}}
\newcommand{\by}{\mathbf{y}}

\newcommand{\bb}{\mathbf{b}}

\newcommand{\BV}{\boldsymbol{V}}

\newcommand{\Bb}{\boldsymbol{b}}

\newcommand{\Bu}{\boldsymbol{u}}
\newcommand{\Bv}{\boldsymbol{v}}

\newcommand{\cP}{\mathcal{P}}
\newcommand{\cS}{\mathcal{S}}

\newcommand{\EUC}{\mathbf E^3}
\newcommand{\SPH}{\bS^3}
\newcommand{\HYP}{\bH^3}

\newcommand{\SOL}{\mathbf{Sol}}

\begin{document}
\pagestyle{myheadings}
\markboth{\centerline{Emil Moln\'ar and Jen\H o Szirmai}}
{On hyperbolic cobweb manifolds}
\title
{On hyperbolic cobweb manifolds
}

\author{Emil Moln\'ar and Jen\H o Szirmai \\
\normalsize Budapest University of Technology and \\
\normalsize Economics Institute of Mathematics, \\
\normalsize Department of Geometry \\
\normalsize Budapest, P. O. Box: 91, H-1521 \\
\normalsize emolnar@math.bme.hu,~szirmai@math.bme.hu
\date{\normalsize{\today}}}

\maketitle
\begin{abstract}
A compact hyperbolic "cobweb" manifold (hyperbolic space form) of symbol
$Cw(6,6,6)$ will be constructed in Fig.1,4,5 as a representant of a presumably
infinite series $Cw(2p,2p,2p)$ $(3 \le p \in \bN$ natural numbers). This is a by-product of our investigations \cite{MSz16}. In that work
dense ball packings and coverings of hyperbolic space $\HYP$ have been constructed on the base of complete hyperbolic Coxeter orthoschemes
$\mathcal{O}=W_{uvw}$ and its extended reflection groups $\bG$ (see diagram in Fig.~3.  and picture of fundamental domain in Fig.~2).
Now $u=v=w=6 (=2p)$.
Thus the maximal ball contained in $Cw(6,6,6)$, moreover its minimal covering bal
l (so diameter) can also be determined. The algorithmic procedure provides us with the proof of our statements.
\end{abstract}
%\tableofcontents

\section{Introduction}

Complete Coxeter simplex (now orthoscheme) $\mathcal{O}=W_{uvw}$, its extended Coxeter's reflection group $\bG$ (diagram in Fig.~3) and its
symmetric Coxeter-Schl\"afli matrix $(i,j\in\{0,1,2,3\}$
\begin{equation}
\begin{gathered}
(b^{ij})=(\langle \bb^i,\bb^j \rangle= (\cos(\pi-\alpha_{ij}))= \\
=\begin{pmatrix}
1&-\cos(\alpha_{01})&-\cos(\alpha_{02})&-\cos(\alpha_{03})  \\
-\cos(\alpha_{10})&1&-\cos(\alpha_{12})&-\cos(\alpha_{13}) \\
-\cos(\alpha_{20})&-\cos(\alpha_{21})&1&-\cos(\alpha_{23}) \\
-\cos(\alpha_{30})&-\cos(\alpha_{31})&-\cos(\alpha_{32})&1
\end{pmatrix}
 \end{gathered} \tag{1.1}
\end{equation}
are mean tools of describing regular polyhedra in absolute geometry,
their metric realization in Euclidean space $\EUC$ (signature of $(b^{ij})$
is $(+,+,+,0))$, in spherical space $\SPH$ (signature is $(+,+,+,+)$), or
in the Bolyai-Lobachevskian hyperbolic space $\HYP$ $(+,+,+,-)$, respectively.

Here projective spherical space  $\cP\cS^3(\bV^4,\BV_4,\bR,\sim^+)$ on real $(\bR)$ vector space $\bV^4$, its dual form space $\BV_4$, model points: $(\bx)=(\by)$ iff $\by=c\bx$, $\bx,\by\in \bV^4$, $c\in\bR^+$ (positive reals), planes: $(\Bu)=(\Bv)$  iff  $\Bv=\frac{1}{c}\Bu$, $\Bu,\Bv \in \BV_4$, $\frac{1}{c}\in \bR^+$. Point $X=(\bx)$ is incident to plane $U=(\Bu)$ iff $\bx\Bu=0$. It stands $X \in u^+$, i.e. $X$ lies in the positive half-space (half-sphere)  $u^+$, iff $\bx\Bu>0$, i.e.  the form $\Bu$ take positive real value on the vector $\bx$. Identifying opposite points $(\bx)$ and $(-\bx)$ (and forms $(\Bu) \sim (-\Bu)$), we get projective space $\cP^3(\bV^4,\BV_4,\bR,\sim)$ for modelling Euclidean and hyperbolic geometry (also elliptic space) in the usual way.

For regular Platonic polyhedron $P(u,v,w)$ we introduce the characteristic Coxeter - Schl\"afli orthoscheme $A_0A_1A_2A_3$ (Fig.2) where $A_3$ is the $3$-centre (or solid centre) of the polyhedron, then $A_2$ is a $2$-centre (face centre) of $P$, $A_1$ is the $1$-centre of an incident (to the previous face) edge, finally $A_0$ is a vertex of the previous edge. At the same time we introduce the side faces $b^0,b^1,b^2,b^3$ of $A_0A_1A_2A_3$, so that $b^i=A_jA_kA_l$, $\{i,j,k,l\}=\{0,1,2,3\}$.  $A_i$-s will be characterized by the vectors $\ba_i\in\bV^4$, while $b^j$-s will be by the forms $\Bb^j\in \BV_4$. The Kronecker symbol
\begin{equation}
\delta_i^j=\ba_i\Bb^j \tag{1.2}
\end{equation}
just describes the incidences of the above simplex as a projective coordinate simplex as well. For the regularity of $P$ we assume (postulate) the angles
\begin{equation}
\begin{gathered}
\angle b^i b^j=\alpha_{ij}, ~\text{so that}~\alpha_{01}=\frac{\pi}{u},~\alpha_{12}=\frac{\pi}{v},~\alpha_{23}=\frac{\pi}{w},\\
(3\le u, v,w \in\bN), ~\text{the others are}~ ~\alpha_{02}=\alpha_{03}=\alpha_{13}=\frac{\pi}{2} ~ \text{(rectangle).} \tag{1.3}
\end{gathered}
\end{equation}
\begin{figure}[ht]

\centering
\includegraphics[width=8cm]{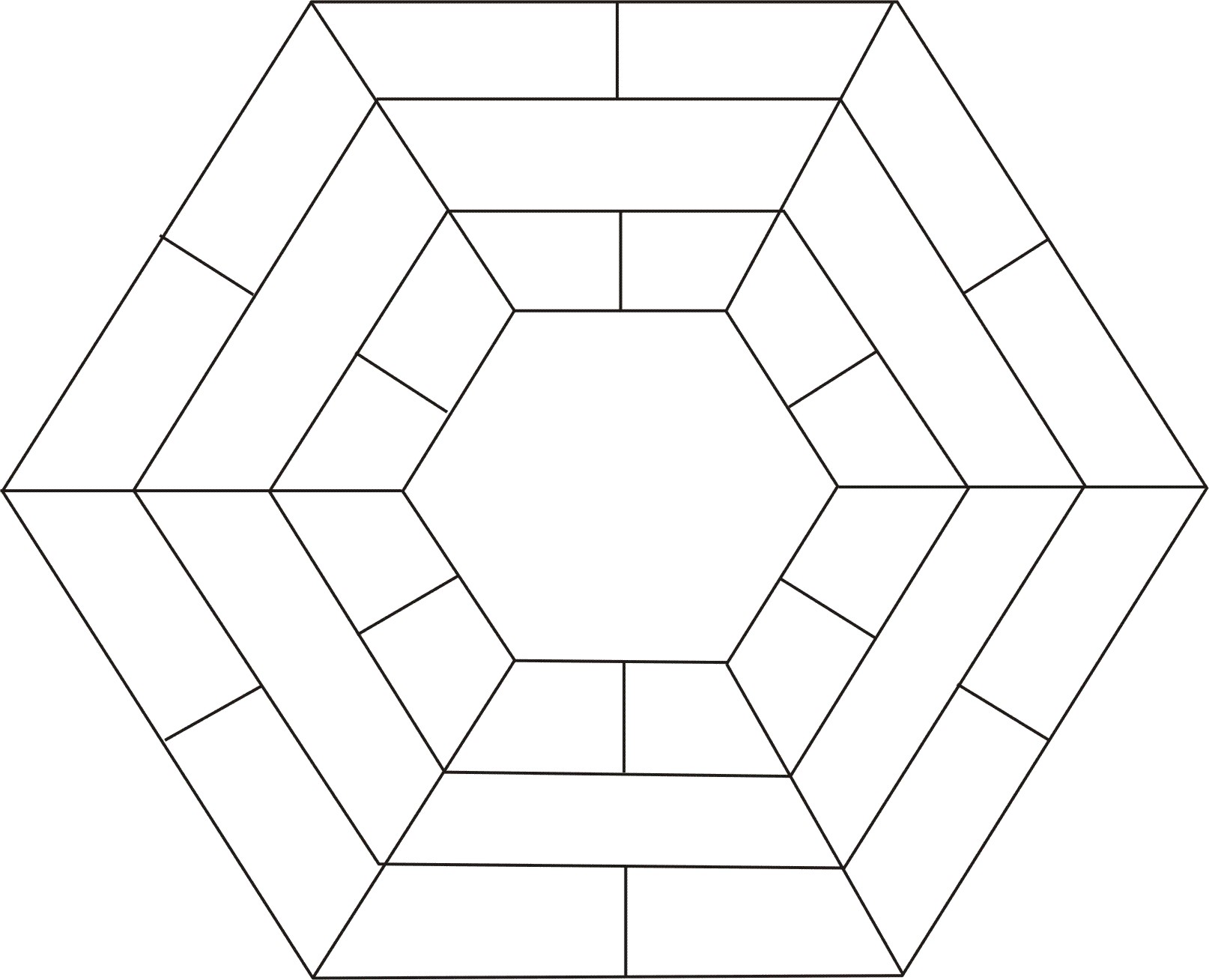}
\caption{Scheme of cobweb polyhedron $Cw(6,6,6)$, $u=v=w=6$, the number of faces is $32$}
\end{figure}

\begin{figure}[ht]
\centering
\includegraphics[width=8cm]{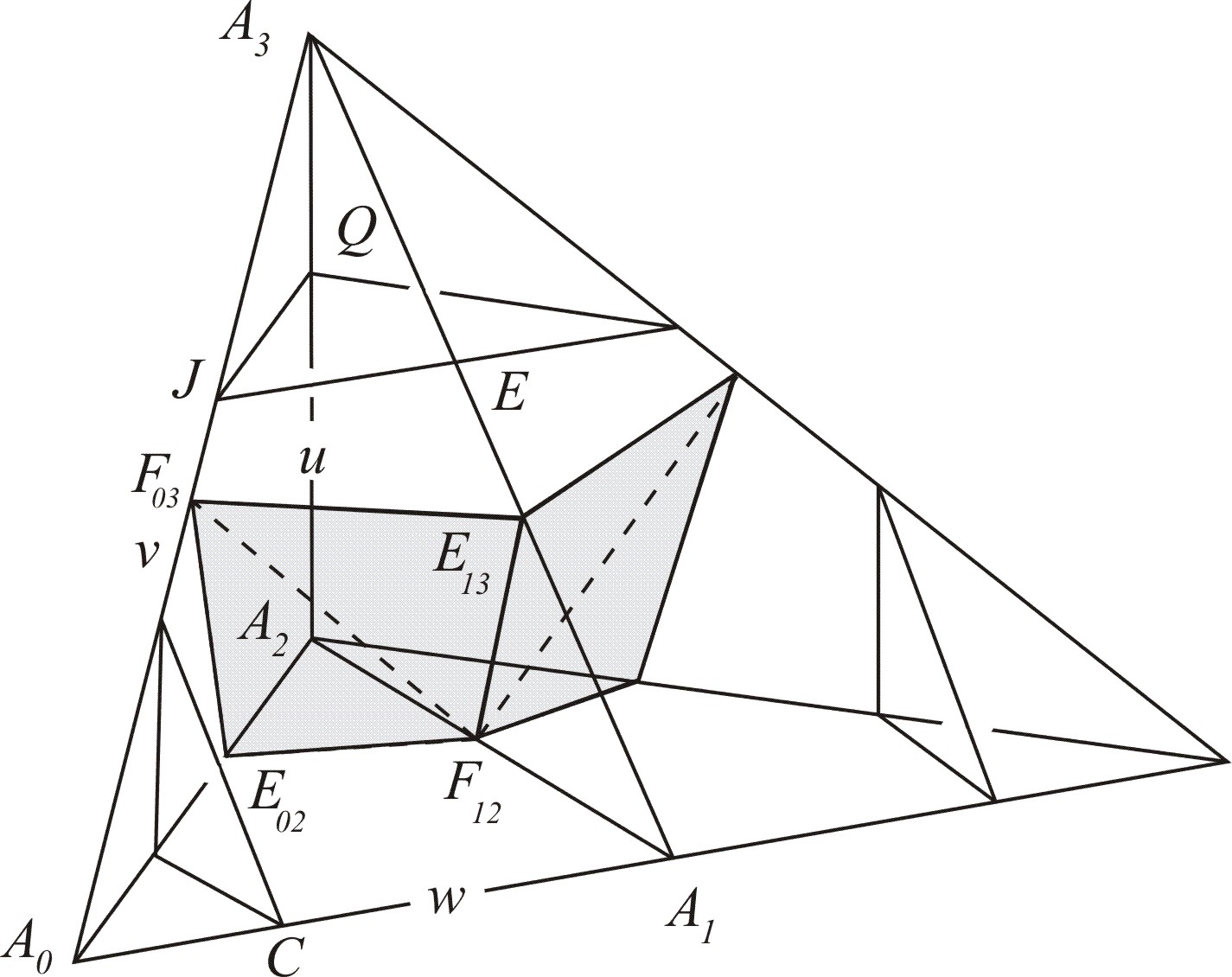}
\caption{Fundamental domains for the half $W_{666}$ and the gluing procedure at point $Q$ for getting $Cw(6,6,6)$}
\end{figure}

Thus, we also guaranteed $A_0A_1A_2A_3=b^0b^1b^2b^3$ to be the characteristic orthoscheme of the regular polyhedron $P(u,v,w)$ with reular $u$-gon faces, meeting $v$ pieces at each vertex. Then
\begin{equation}
\frac{\pi}{u}+\frac{\pi}{2}+\frac{\pi}{v}>\pi \Leftrightarrow \frac{\pi}{u}+\frac{\pi}{v}>\frac{\pi}{2}, \tag{1.4}
\end{equation}
for the angle sum of spherical triangle, guarantees $A_3$ as proper centre. This is equivalent with
\begin{equation}
\det\left[ \begin{array}{ccc}
1 &-\cos\frac{\pi}{u} & 0 \\
-\cos \frac{\pi}{u} &1& -\cos\frac{\pi}{v} \\
0 & -\cos\frac{\pi}{v} &1
\end{array} \right]>0 ~\text{and}~(+,+,+,\cdot) \tag{1.5}
\end{equation}
in the signature of the very important quadratic form
\begin{equation}
\begin{gathered}
b^{ij}u_iu_j=(b^{00}u_0u_0+2b^{01}u_0u_1+b^{11}u_1u_1+2b^{12}u_1u_2+b^{22}u_2 u_2)+\\
+2b^{23}u_2u_3+b^{33}u_3u_3. \tag{1.6}
\end{gathered}
\end{equation}
(Einstein-Schouten convention is and will be applied).

Similarly, the vertex $A_0$ will be proper iff
\begin{equation}
\frac{\pi}{v}+\frac{\pi}{w}>\frac{\pi}{2}. \tag{1.7}
\end{equation}
Thus, the face angle $\alpha_{23}=\frac{\pi}{w}$ will be important for $P$, so that it will be a space-filler regular polyhedron, $w$ pieces meeting at each edges.

Now the hyperbolic regular mosaic (honeycomb) with congruent pieces of $P$ will be characterized just the determinant of (1.1) specialized
\begin{equation}
\det\left[ \begin{array}{cccc}
1 &-\cos\frac{\pi}{u} & 0&0 \\
-\cos \frac{\pi}{u} &1& -\cos\frac{\pi}{v}&0 \\
0 & -\cos\frac{\pi}{v} &1&-\cos \frac{\pi}{w} \\
0&0&-\cos \frac{\pi}{w} &1
\end{array} \right]<0. \tag{1.8}
\end{equation}

For instance, the famous space-filler dodecahedron $P(5,3,5)$, with $5$ pieces at each edge can be described in this manner as in \cite{WeSe33}. Then the $(+,+,+,-)$ signature of the quadratic form (1.6) guarantees $\HYP$ as relevant space in projective-metric realization (Beltrami-Cayley-Klein model).

Now we consider doubly truncated orthoschemes, as in Fig.2, where
\begin{equation}
\frac{\pi}{u}+\frac{\pi}{v}<\frac{\pi}{2}, ~ \text{change also the sign in (1.5)} \tag{1.4'}
\end{equation}
guarantees that $A_3$ is outer point. Then we cut (truncate) the simplex with its polar plane $a_3=QEJ$. Furthermore, we assume that
\begin{equation}
\frac{\pi}{v}+\frac{\pi}{w}<\frac{\pi}{2},  \tag{1.7'}
\end{equation}
i.e. $A_0$ is also outer vertex with truncating polar plane $a_0$. The corresponding minor determinant is also negative, as the complete determinant in (1.8), too.

\begin{figure}[ht]
\centering
\includegraphics[width=7cm]{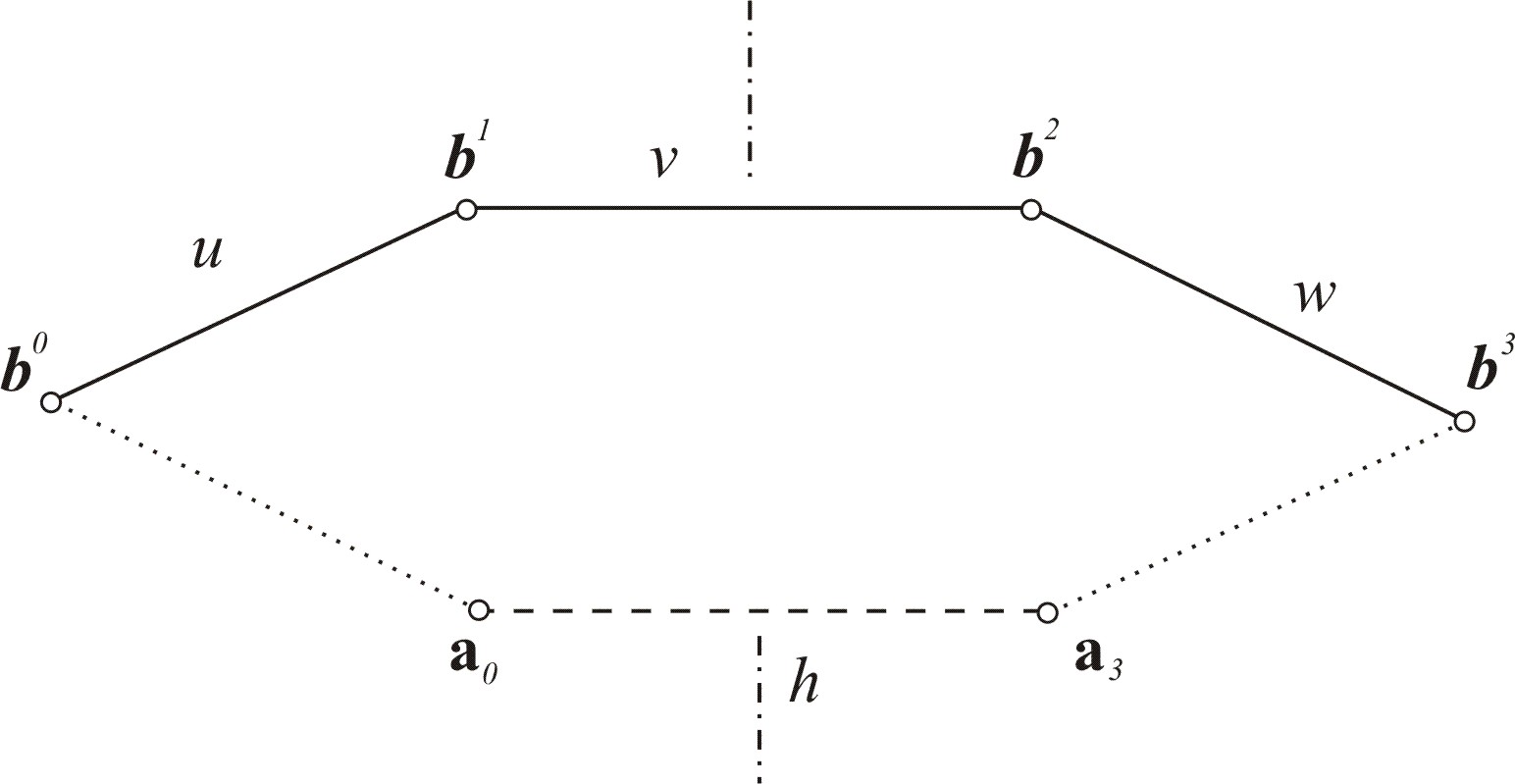}
\caption{The general Coxeter-Schl\"afli diagram for the (half-turn) extended complete Coxeter orthoscheme group $\bG(u=w,v)$}
\end{figure}
Our aim in \cite{M84} was to construct the first non-orientable series (besides orientable ones) on generalized regular polyhedron
$P_Q$ with centre $Q$ in Fig.~2, where $4u$ truncated orthoschemes meet. Moreover, we intended to equip this $P_Q$ with face identification (gluing procedure) so that this $\tilde{P}_Q$ has hyperbolic ball-like neighborhood in each point. To that procedure $u=v=w$ even was needed.
{\it But then the complete (truncated) orthoscheme $\mathcal{O}=W_{uuu}$ have
an additional halfturn $\mathbf{h}$ about axis $F_{03}F_{12}$, where $F_{03}$ and
$F_{12}$ are midpoints of $A_0A_3$ and $A_1A_2$, respectively.
So the problem arises, whether a smaller hyperbolic space form (a half of
the previous one) can also be constructed. Our results are shortly formulated in}
\begin{theorem}
The cobweb manifold $Cw(6,6,6)$ to Fig.1,2 has been
constructed by face identification in Fig. 4,5.

The fundamental group $\mathbf{Cw}(6,6,6)$ can be described by
$3$-generators and two relations in formulas (2.14-15).

The volume of $Cw(6,6,6)$ is $\approx 8.29565$ in (3.7). The largest ball contained in $Cw(6,6,6)$ is of radius
$r \approx 0.57941$. The diameter of $Cw(6,6,6)$ is $2R \approx 3.67268$ by (3.4-5).
\end{theorem}
The usual diagram in Fig.3 describes the extended complete orthoscheme group $\bG$, where $\Bb^i$, for plane reflections $\mathbf{m}^i$, and $\ba_0$,  $\ba_3$ are indicated by $6$ nodes. The branches $3 \le u,v,w \in \bN$ indicate the order of reflection products (as rotations). No branch means orthogonality of reflection planes, i.e. order $2$ of rotation. The truncating planes $a_0$, $a_3$ and half turn $\mathbf{h}$ are also indicated in the usual way. Dotted and dashed lines sign planes with common perpendiculars.

The half truncated orthoscheme $W_{uvw}$ is also indicated in Fig.~2, where the gluing procedure at $Q=A_2A_3\cup a_3$ are illustrated, too. After having glued $4u=24$ half truncated orthoscheme at $Q$ we get cobweb as edge frame of the new polyhedron $P_Q$ as the fundamental domain for the space form $Cw(6,6,6)$. We found similar hyperbolic space forms $Cw(10,10,10)$ and $Cw(14,14,14)$, but not yet for $u=v=w=8$.

{\it It seems to be an interesting open problem to give a unified
construction scheme for all conjectured cobweb manifolds $Cw(2p,2p,2p)$
($3<p\in \bN$).}
\begin{figure}[ht]
\centering
\includegraphics[width=8cm]{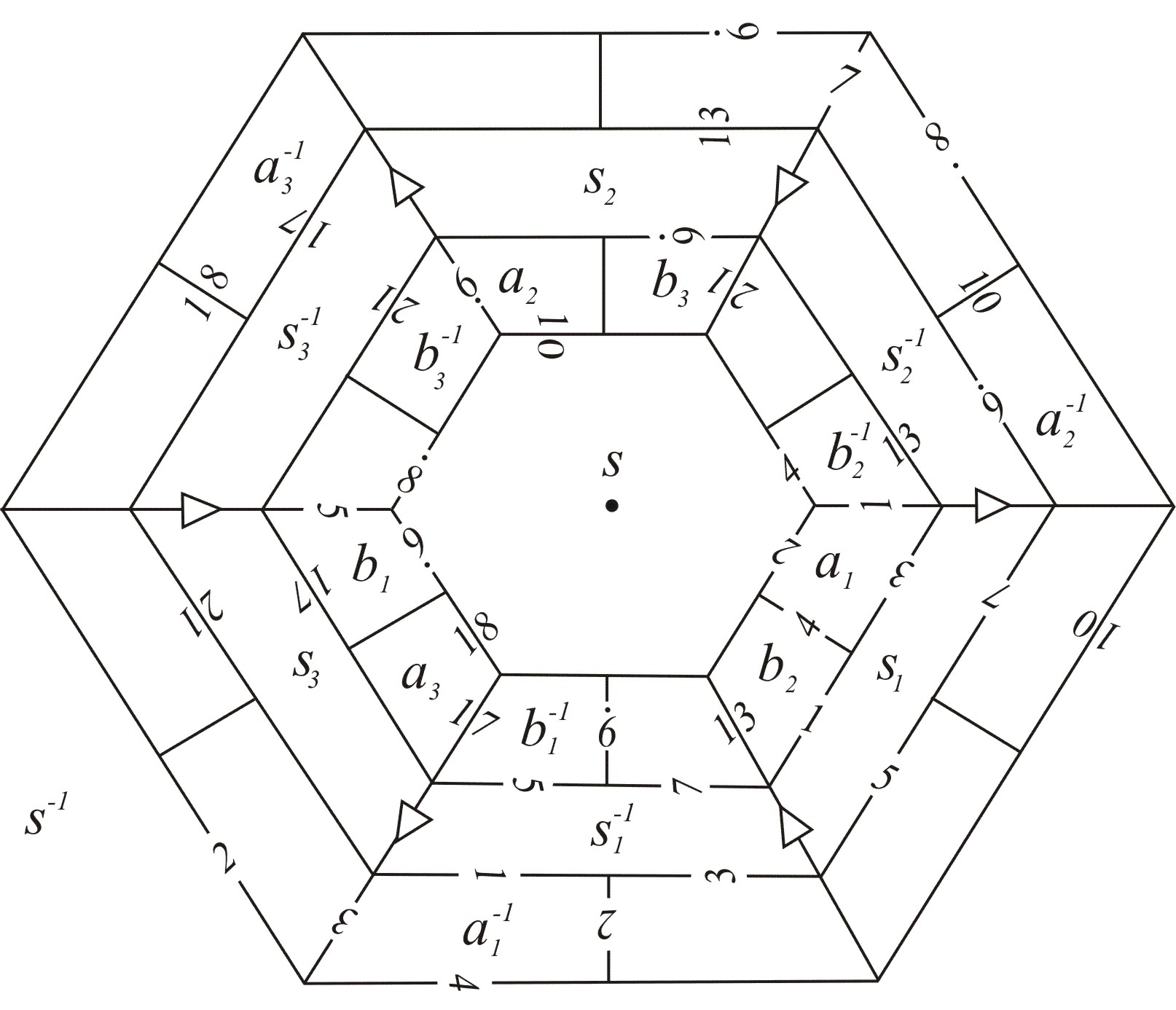}
\caption{The strategy for face pairing motions to the cobweb polyhedron $Cw(6,6,6)$}
\end{figure}

We give a short list of references only, where the interested Reader can find further items. After the classical work \cite{WeSe33}, the first hyperbolic space form series were constructed in \cite{Be71}, \cite{AJ80}, \cite{Gu80}, \cite{Ve84}. We advise to study \cite{M12} for a more popular introduction and applications. In \cite{CaTe09} you find also a graph interpretation of "football manifolds" of the first author. \cite{CaMoSpSz} is an example of space forms in other Thurston geometries, namely in $\SOL$. \cite{P} discusses algorithmic problems for polyhedra to represent space forms. \cite{MSz} describes these $8$ homogeneous $3$-spaces in a popular projective interpretation. \cite{W} is a classical book of space form problem.
\section{Construction of cobweb manifold $Cw(6,6,6)$}
By the theory of \cite{W} we have to construct a fixed point free group acting in hyperbolic space $\HYP$ with compact fundamental domain.
\begin{figure}[ht]
\centering
\includegraphics[width=10cm]{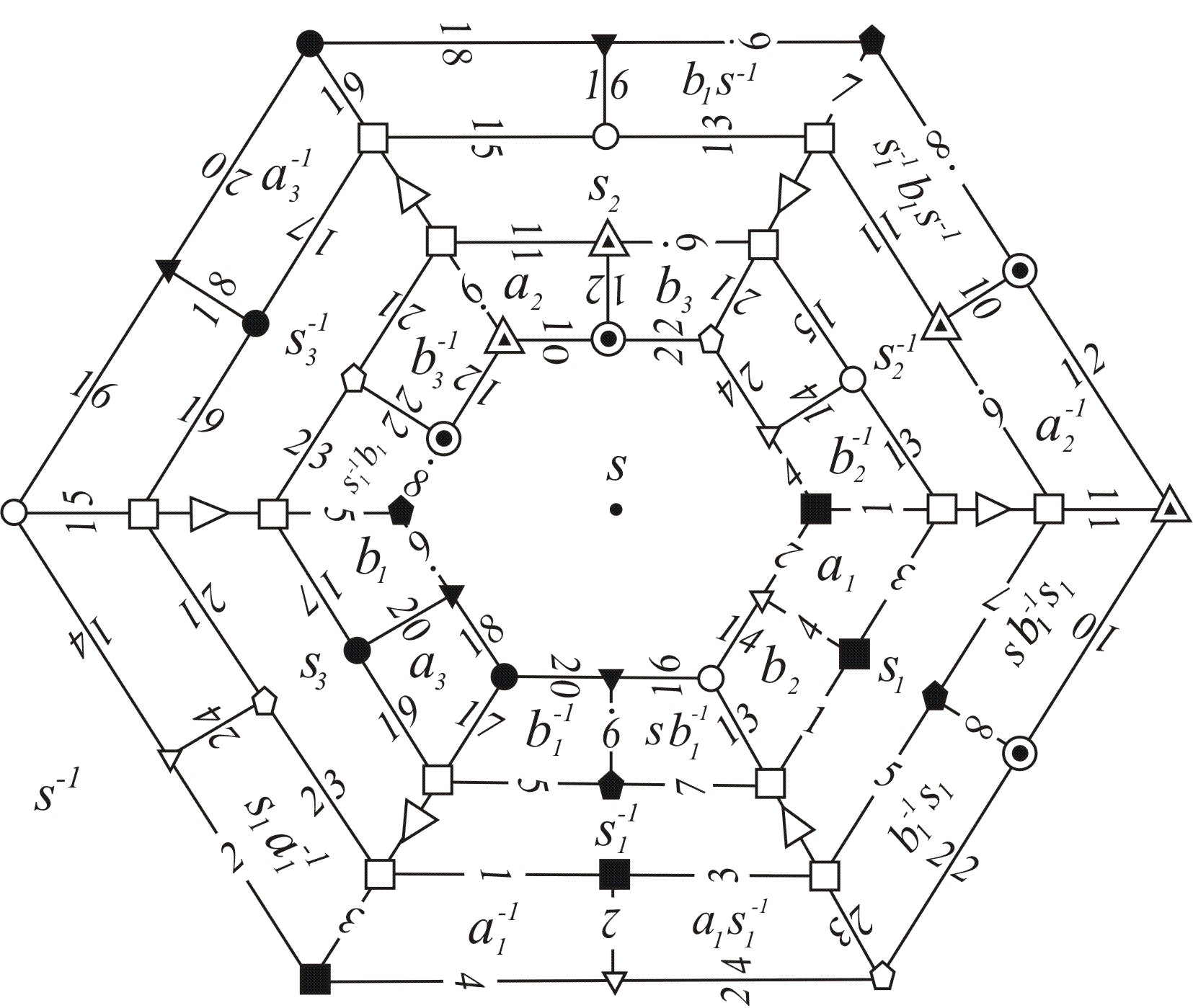}
\caption{The hyperbolic cobweb manifold $Cw(6,6,6)$ with complete (symbolic) face pairing.}
\end{figure}
In the Introduction (Sect. 1) to Fig.~2,3 we have described the extended reflection group $\bG(6,6,6)$ to the complete Coxeter orthoscheme $\mathcal{O}=W_{666}$ and glued together the cobweb polyhedron $Cw(6,6,6)$ as Dirichlet-Voronoi (i.e. D-V) cell of the kernel point $Q$ by its orbit under the group $\bG(6,6,6)$ by Fig.~3. Now by Fig.~1,4,5 we shall give the face identification of $Cw(6,6,6)$, so that it will be fundamental polyhedron of the fixed-point-free group, denoted also by $\mathbf{Cw}$, generated just by the face identifying isometries (as hyperbolic screw motions).

The construction scheme will be illustrated in Fig.4,5. The strategy has been indicated in Fig.4, while in Fig.5 the complete $\tilde{C}w(6,6,6)$ has appeared with face pairs, signed and numbered (from $1$ to $24$) edge triples, signed vertex classes, all together $1+3 \times 3=10$ vertex classes.

In the previous considerations (Sect.1 ) we described the extended complete reflection group $\bG(6,6,6)$. Fig.3, by its fundamental domain, where the stabilizer subgroups of $\bG$ can also be established. By gluing $4u=24$ domains at $Q$ (whose stabilizer subgroup $\bG_Q$ is just of order $|stab_Q\bG|=4u=24$) we can "kill out" the fixed points of $\bG$.

To this $v=u=6$ edge domains (signed by $-\!\!-\!\!\!\vartriangleright$) is just necessary and sufficient at the former edge $F_{03}J$ of  half $W_{666}$ for a ball-like neighbourhood of points in $-\!\!-\!\!\!\vartriangleright$ edges. This can be achieved by $3$ screw motions $\bs_1,\bs_2,\bs_3$ for the $6$ middle faces of the cobweb $Cw(6,6,6)$, $\bs_i:s_i^{-1} \rightarrow s_i$. The $12$ images of the former $F_{03}$ will form a vertex class $\square$, since $|stab_{F_{03}}\bG|=4v=24$ and just $24(=4u)$ domains will form the ball-like neighbourhood at these $12$ $\square$-images.

The most crucial roles are played by the former edges at the halving planes of the orthoscheme $W_{666}$ to the half-turn axis $F_{03}F_{12}$. The stabilizers of the mirror points are of order $2$ divided into two parts at $W_{666}$, namely at $F_{03}E_{02}$ and at $F_{03}E_{13}$ for the odd numbered edges $1,3,\dots 21,23$, and of $F_{12}E_{02}$, $F_{12}E_{13}$ for the even numbered edges $2,4,\dots,22,24$, respectively. The different roles of reflection mirrors of $m^1$ and $m^2$, resp. $m^0$ and $m^3$ in the gluing procedure at $Q$ yield that both edge classes appear in three copies on the cobweb polyhedron $Cw(6,6,6)$, each of both classes maintains ball-like neighborhood at each point of them.

Now comes our tricky constructions for identifying the former  half-turn faces, furthermore  the two base faces $s^{-1}$ and $s$ of $Cw(6,6,6)$ with each other (see Fig.2 and Fig.4,5). Two from the edge triple $1$ (to $F_{03}E_{13}$) lie on the faces $s_1^{-1}$ and $s_1$, we introduce the deciding third edge $1$ (to $F_{03}E_{12}$) and the orientation preserving motions $\ba_1:a_1^{-1}\rightarrow a_1$ and its inverse $\ba_1^{-1}:a_1\rightarrow a_1^{-1}$ by the mapping faces $a_1^{-1}$, $a_1$. This edge triple $1$ with faces $s_1^{-1},s_1,a_1^{-1},a_1$ defines a third face pairing identification $\bb_2:b_2^{-1}\rightarrow b_2$ so that
\begin{equation}
1:~\ba_1^{-1}\bs^1=\bb_2 ~ \text{holds.} \tag{2.1}
\end{equation}
Namely, three image polyhedra join each other, e.g. at the first $1$ edge in space $\HYP$ (now in combinatorial sense):
\begin{equation}
a_1^{-1}(Cw)s_1^{-1},~s_1^{\bs_1^{-1}}(Cw) ^{\bs_1^{-1}}b_2^{\bs_1^{-1}},~(b_2^{-1})^{\bb_2\bs_1^{-1}}(Cw) ^{\bb_2\bs_1^{-1}}a_1^{\bb_2\bs_1^{-1}}. \tag{2.2}
\end{equation}
Now comes again the identity polyhedron through the images
\begin{equation}
(a_1^{-1})^{\ba_1\bb_2\bs_1^{-1}}(Cw) ^{\ba_1\bb_2\bs_1^{-1}}(s_1^{-1})^{\ba_1\bb_2\bs_1^{-1}}. \tag{2.3}
\end{equation}
That means we get
\begin{equation}
\ba_1\bb_2\bs_1^{-1}=\mathbf{1} ~ \text{the identity}~ \Leftrightarrow \bb_2=\ba_1^{-1}\bs_1, \tag{2.4}
\end{equation}
indeed as in (2.1).

This general method for space filling with fundamental polyhedron, equipped by face pairing generated group, has been discussed in \cite{M92} in more details.

The procedure is straightforward in Fig.4, then in Fig.5. The next edge class $3$ defines the new face pairs with motion
\begin{equation}
3:~\bs_1\ba_1^{-1}:~ a_1s_1^{-1} \rightarrow s_1a_1^{-1}  ~ \text{and its inverse}~ \ba_1\bs_1^{-1}:~ s_1a_1^{-1} \rightarrow a_1s_1^{-1}. \tag{2.5}
\end{equation}
The next edge class $2$ (to edges $F_{12}E_{13}$, $F_{12}E_{02}$ in Fig.2), again with triple, just defines the identification of the base faces of $Cw(6,6,6)$:
\begin{equation}
2:~\bs=(\ba_1\bs_1^{-1})\ba_1:~ s^{-1} \rightarrow s.  \tag{2.6}
\end{equation}
For the edge triple $4$ (to $F_{12}E_{13}$, $F_{12}E_{02}$ in Fig.2) we get relation
\begin{equation}
\begin{gathered}
4:~\ba_1\bb_2^{-1}(\bs_1^{-1})=\mathbf{1}~ \text{which will be}~ \ba_1(\bs_1^{-1}\ba_1)(\ba_1^{-1}\bs_1\ba_1^{-1})=1\\
\text{a trivial relation for generators}~ \ba_1 ~\text{and}~ \bs_1. \tag{2.7}
\end{gathered}
\end{equation}
Our next "lucky" choice (it comes from the starting triples $1$, see at the end of Sect. 2) for edge triple $5$ was influenced by the trigonal symmetry of our cobweb polyhedron $Cw(6,6,6)$. Then the triples $6,7,8$ follow as formulas in (2.8) show
\begin{equation}
5:~\bs_1^{-1}\bb_1,~6:~\bb_1\bs^{-1},~7:~(\bs \bb_1^{-1})\bs_1,  ~ 8:~ ( \bs_1^{-1}\bb_1)\bs^{-1}(\bs\bb_1^{-1}\bs_1)=\mathbf{1}. \tag{2.8}
\end{equation}
The last relation to edge class $8$ is trivial again.

The procedure is straightforward now, and it nicely closes. The edge triple $9$ defines the face pairing motion
\begin{equation}
\ba_2:~ a_2^{-1} \rightarrow a_2  ~ \text{moreover, a new motion}~ \ba_2^{-1}\bs_2=\bb_3.  \tag{2.9}
\end{equation}
The further triples $10-16$ and identifications are completely analogous. The same holds for edge triples $17-24$, starting with the face pairing motion $\ba_3$.

It turns out that the first three screw motions $\bs_1, \bs_2, \bs_3$  can be expressed by $\ba_1, \ba_2, \ba_3$  and $\bs$ at triples $2, 10, 18$, respectively:
\begin{equation}
\bs_1=\ba_1\bs^{-1}\ba_1,~ \bs_2=\ba_2 \bs^{-1}\ba_2,~\bs_3=\ba_3\bs^{-1}\ba_3.  \tag{2.10}
\end{equation}
The relation, to the middle edge class $-\!\!\!-\!\!\!\vartriangleright$ of $6$ edges yields then the relation
\begin{equation}
\begin{gathered}
\mathbf{1}=(\ba_1\bs^{-1}\ba_1)^2(\ba_2\bs^{-1}\ba_2)^2(\ba_3\bs^{-1}\ba_3)^2
\end{gathered} \tag{2.11}
\end{equation}
for the fundamental group of our cobweb manifold $Cw(6, 6, 6)$.

But in this cyclic process, the pairing motion $\bb_1$ to edge class $5$ is not independent. Symilarly to $1:~\bb_2 =\ba_1^{-1}\bs_1=\bs^{-1}\ba_1$, as above, we cyclically obtain $\bb_1=\ba_3^{-1}\bs_3=\bs^{-1}\ba_3$. So we get, at the edge class $7$, the motion $(\bs \bb_1^{-1})\bs_1=\bs(\ba_3^{-1}\bs)(\ba_1\bs^{-1}\ba_1)$,
\begin{equation}
\begin{gathered}
\mathbf{1}=(\bs\ba_3^{-1}\bs\ba_1\bs^{-1}\ba_1)\bs\ba_2^{-1},  \tag{2.12}
\end{gathered}
\end{equation}
as well at (2.11). Analogously, at edge class $19$, so we get the relations
\begin{equation}
\begin{gathered}
\mathbf{1} = (\bs\ba_1^{-1}\bs\ba_2\bs^{-1}\ba_2)\bs\ba_3^{-1};~\text{and} \\
\mathbf{1} = \bs \ba_2^{-1} \bs \ba_3 \bs{-1} \ba_3 \bs \ba_1^{-1}; ~ \text{eliminating}~ \ba_3, ~ \text{we get} \\
\mathbf{1} = \bs a_2^{-1} \bs (\bs \ba_1 \bs^{-1} \ba_1 \bs \ba_2^{-1} \bs) \bs^{-1} (\bs \ba_1^{-1} \bs \ba_2 \bs^{-1} \ba_2 \bs)\bs \ba_1^{-1}, \\
\text{so} ~ \mathbf{1}_{18} = \bs \ba_2^{-1} \bs_2 \ba_1 \bs^{-1} \ba_1 \bs \ba_2^{-1} \bs \ba_1^{-1} \bs \ba_2 \bs^{-1} \ba_2 \bs_2 \ba_1^{-1}.
\end{gathered} \tag{2.13}
\end{equation}
Thus,$\ba_3$ can be expressed, from (2.12) and (2.13), keeping $1\leftrightarrow 2$ symmetry, and we obtain for the fundamental group $\mathbf{Cw}(6, 6, 6)$ three relations for the three generators $\ba_1$, $\ba_2$, $\bs$. From (2.12-13) we get first (a $10$ letters relation)
\begin{equation}
\begin{gathered}
\mathbf{1}_{10}=\ba_1\ba_1\bs^{-1}\ba_1\bs\ba_2^{-1}\ba_2^{-1}\bs\ba_2^{-1}\bs^{-1},
\end{gathered} \tag{2.14}
\end{equation}
         
Second, from (2.11) we obtain (a $38$ letters relation), symmetrically with indices $1$ and $2$,
\begin{equation}
\begin{gathered}
\mathbf{1}_{38}=(\ba_1\bs^{-1}\ba_1)^2(\ba_2\bs^{-1}\ba_2)^{2}
(\bs\ba_1\bs^{-1}\ba_1\bs\ba_2^{-1}\bs\ba_1^{-1}\bs\ba_2\bs^{-1}\ba_2\bs)^2.
\end{gathered} \tag{2.15}
\end{equation}
But we do not give more details, e.g the first homology group $\mathbf{H}_1$ of manifold ${Cw}(6, 6, 6)$, by the commutator factorgroup of $\mathbf{Cw}(6, 6, 6)$ can easily be determined, as $\mathbf{Z}_3 \times \mathbf{Z}_{12} \times \mathbf{Z}_6$ direct product of cyclic groups.

Of course, this group $\mathbf{Cw}$, is a subgroup of our former $\bG(6, 6, 6)$ by Fig.3. These generators $\ba_1$, $\ba_2$, $\bs$  can be expressed by the former reflections $\mathbf{m}^0$, $\mathbf{m}^1$, $\mathbf{a}_3$ and the half-turn $\mathbf{h}$
about $F_{03}F_{12}$ (Fig.2,3).
\section{Some metric properties of our cobweb manifold $Cw(6,6,6)$}
Hyperbolic metrics, defined by the Coxeter-Schl\"afli matrix (1.1), now is specialized by $u = v = w = 6$ in (1.3) and Fig. 3. As it is well-known, the Kronecker relations (1.2) involve that the inverse matrix of (1.1) for the coordinate orthoscheme $A_0 A_1 A_2 A_3 = b^0 b^1 b^2 b^3$, now in the form of (1.8), determines the distance metrics of $\HYP$ just by the scalar product (of signature $(+,+,+,-)$)
\begin{equation}
\begin{gathered}
(a_{ij})=(b^{ij})^{-1}=\langle \ba_i, \ba_j \rangle:=\\
=\frac{1}{B} \begin{pmatrix}
\sin^2{\frac{\pi}{w}}-\cos^2{\frac{\pi}{v}}& \cos{\frac{\pi}{u}}\sin^2{\frac{\pi}{w}}& \cos{\frac{\pi}{u}}\cos{\frac{\pi}{v}} & \cos{\frac{\pi}{u}}\cos{\frac{\pi}{v}}\cos{\frac{\pi}{w}} \\
\cos{\frac{\pi}{u}}\sin^2{\frac{\pi}{w}} & \sin^2{\frac{\pi}{w}} & \cos{\frac{\pi}{v}}& \cos{\frac{\pi}{w}}\cos{\frac{\pi}{v}} \\
\cos{\frac{\pi}{u}}\cos{\frac{\pi}{v}} & \cos{\frac{\pi}{v}} & \sin^2{\frac{\pi}{u}}  & \cos{\frac{\pi}{w}}\sin^2{\frac{\pi}{u}}  \\
\cos{\frac{\pi}{u}}\cos{\frac{\pi}{v}}\cos{\frac{\pi}{w}}  & \cos{\frac{\pi}{w}}\cos{\frac{\pi}{v}} & \cos{\frac{\pi}{w}}\sin^2{\frac{\pi}{u}}  & \sin^2{\frac{\pi}{u}}-\cos^2{\frac{\pi}{v}}
\end{pmatrix},
\end{gathered} \tag{3.1}
\end{equation}
where
$$
B=\det(b^{ij})=\sin^2{\frac{\pi}{u}}\sin^2{\frac{\pi}{w}}-\cos^2{\frac{\pi}{v}} <0, \ \ \text{i.e.} \ \sin{\frac{\pi}{u}}\sin{\frac{\pi}{w}}-\cos{\frac{\pi}{v}}<0.
$$
in formula (1.8) is crucial. Thus the (complex) distances $A_iA_j$ of orthoscheme vertices $A_i(\ba_i)$, $A_j(\ba_j)$ can be calculated by
\begin{equation}
\cosh\Big(\frac{A_iA_k}{k}\Big)=\frac{-a_{ij}}{\sqrt{a_{ii}a_{jj}}} \tag{3.2}
\end{equation}
e.g. for the proper vertices $A_1$,  $A_2$
\begin{equation}
\cosh\Big(\frac{A_1A_2}{k}\Big)=\frac{\cos\Big(\frac{\pi}{v}\Big)}{{\sqrt{-B} \sin\Big(\frac{\pi}{w}\Big)\sin\Big(\frac{\pi}{u}\Big)}} =2\sqrt{3} \rightarrow A_1A_2 \approx 1.91408 \notag
\end{equation}
where  $a_{11}$, $a_{22}$, $a_{12}$ < 0, then the other distances.
The points of the truncating polar plane $a_3$,  to the outer vertex $A_3(\ba_3)$, e.g. $Q(\mathbf{q})$ is characterized by the conjugacy
\begin{equation}
\begin{gathered}
\langle \mathbf{q}, \ba_3 \rangle=0,~\text{so}~ \mathbf{q}\sim a_{33}\ba_2-a_{23} \ba_3~ \text{and} ~
\cosh\Big(\frac{QA_2}{k}\Big)=\frac{1}{\sqrt{a_{33}}}
\end{gathered} \tag{3.3}
\end{equation}
can also be calculated. The natural length unit $k =\sqrt{-1/K}$ ($K$  is the constant negative sectional curvature) can be chosen to $k = 1$. So $2\times QA_2$ is the height of our polyhedron $Cw(6, 6, 6)$.

In paper \cite{MSz16} we systematically computed (by computer program) the radius $r$ of the maximal packing ball (inscribed into $Cw(6, 6, 6)$. Namely,
\begin{equation}
\begin{gathered}
r=\min\{QA_2,Qb^2=QE,Q(F_{03}F_{12})\} \\
r=QA_2=\mathrm{arcosh}\Big(\frac{1}{\sqrt{a_{33}}}\Big) =\\
=\mathrm{arcosh}\Big(\sqrt{1+\frac{\sin^2\frac{\pi}{u}\cos^2\frac{\pi}{w}}{\cos^2\frac{\pi}{v}-\sin^2\frac{\pi}{u}}}\Big) \approx 0.57941.
\end{gathered} \tag{3.4}
\end{equation}
This serves the maximal ball, contained in the manifold.
Similarly, the covering ball of $Cw(6, 6, 6)$ can be determined by the radius
\begin{equation}
\begin{gathered}
R=\max\{QF_{12},QF_{03},QE_{13}\} \\
R=QF_{03}=\mathrm{arcosh}\Big( \frac{\cot\frac{\pi}{u}\cos\frac{\pi}{v}}{\sqrt{2Ba_{33}(a_{33}+a_{03})}}\Big) \approx 1.83634.
\end{gathered} \tag{3.5}
\end{equation}
$2R \approx 3.67268$ is the diameter of our manifold.
The volume of the complete truncated orthoscheme $W_{666}$ can be computed by a general formula of R. Kellerhals \cite{K89} by the ideas of N.I. Lobachevsky:
\begin{theorem}{\rm{(R.~Kellerhals)}} The volume of a three-dimensional hyperbolic
complete orthoscheme $\mathcal{O}=W_{uvw} \subset \mathbf{H}^3$
is expressed with the essential
angles $\alpha_{01}=\frac{\pi}{u}$, $\alpha_{12}=\frac{\pi}{v}$, $\alpha_{23}=\frac{\pi}{w}$, $(0 \le \alpha_{ij}
\le \frac{\pi}{2})$
(See our formulas (1.1), (1.3), (1.8), (3.1) and Fig.2) in the following form:

\begin{align}
&\mathrm{Vol}(\mathcal{O})=\frac{1}{4} \{ \mathcal{L}(\alpha_{01}+\theta)-
\mathcal{L}(\alpha_{01}-\theta)+\mathcal{L}(\frac{\pi}{2}+\alpha_{12}-\theta)+ \notag \\
&+\mathcal{L}(\frac{\pi}{2}-\alpha_{12}-\theta)+\mathcal{L}(\alpha_{23}+\theta)-
\mathcal{L}(\alpha_{23}-\theta)+2\mathcal{L}(\frac{\pi}{2}-\theta) \}, \tag{3.6}
\end{align}
where $\theta \in [0,\frac{\pi}{2})$ is defined by:
$$
\tan(\theta)=\frac{\sqrt{ \cos^2{\alpha_{12}}-\sin^2{\alpha_{01}} \sin^2{\alpha_{23}
}}} {\cos{\alpha_{01}}\cos{\alpha_{23}}},
$$
and where $\mathcal{L}(x):=-\int\limits_0^x \log \vert {2\sin{t}} \vert dt$ \ denotes the
Lobachevsky function.
\end{theorem}
The $4\times 6/2 = 12$ times of $\mathrm{Vol}(W_{666})$ is equal to the volume of our cobweb manifold
\begin{equation}
\mathrm{Vol}(Cw(6, 6, 6)) = 12 W_{666} \approx 8.29565.  \tag{3.7}
\end{equation}
The densities of packing and covering are
\begin{equation}
\begin{gathered}
\delta(6, 6, 6) = \mathrm{Vol}(B(r))/\mathrm{Vol}(Cw) = 0.10503 , \\ \Delta(6, 6, 6) = \mathrm{Vol}(B(R))/\mathrm{Vol}(Cw) = 6.05670, \end{gathered} \tag{3.8}
\end{equation}
respectively, play relevant roles for our manifold.

Computer programs of the second author can provide further fruitful results.

\end{document}